\documentclass[12pt]{svmult}
\usepackage[french,english]{babel}
\usepackage[T1]{fontenc}
\usepackage[applemac]{inputenc}
\usepackage{amssymb,url,xspace,amsmath,eucal,yhmath}
\usepackage{mathrsfs,eufrak,array,epsfig,psfrag,graphicx,dsfont,bbm,stmaryrd}
\usepackage[all]{xy}
\usepackage{hyperref}
\usepackage{graphicx}

\renewcommand\Pr[1]{\mathbb{P}\left[#1\right]}

\newcommand\Prmu[1]{\mathbb{P}_{\mu}\left[#1\right]}

\newcommand\Es[1]{\mathbb{E}\left[#1\right]}
\newcommand\Esmu[1]{\mathbb{E}_{\mu}\left[#1\right]}

\newcommand  \Esmuzgeqn[1] {\mathbb{E}_\mu\left[ #1 \, | \, \zt\geq n\right]}

\newcommand \fl[1] {\left\lfloor #1 \right\rfloor}

\def \Card {\mathop{{\rm Card } }\nolimits}

\def \P {\mathbb{P}}

\def \Pmu {\mathbb{P}_\mu}
\def \Pmuzn {\mathbb{P}_\mu[ \, \cdot \, | \, \zt=n]}

\def \Pmuzgeqn {\mathbb{P}_\mu[ \, \cdot \, | \, \zt\geq n]}

\def \N {\mathbb N}
\def \E {\mathbb E}
\def \D {\mathbb D}

\def \C {\mathcal{C}}
\def \R {\mathbb R}

\def \D {\mathbb D}
\def \T {\mathbb T}
\def \Z {\mathbb Z}
\def \F {\mathcal{F}}

\def \mE {\mathcal{E}}

\def \W {\mathcal{W}}

\def \l {\lambda}

\def \H {H^\textrm{exc}}

\def \d {\displaystyle}
\def \a {\alpha}
\def \b {\beta}

\def \e {\epsilon}

\def \Nn {\textnormal{\textbf{N}}}

\def \GW {\textnormal{GW}}
\def \z {\zeta}

\def \t {\mathfrak{t}}

\def \| { \, | \,}

\def \zt {\zeta(\tau)}

\renewcommand\i[1]{\inf_{\left[#1\right]}}

\long\def\symbolfootnote[#1]#2{\begingroup%
\def\thefootnote{\fnsymbol{footnote}}\footnote[#1]{#2}\endgroup}

\renewcommand{\keywords}[1]{ \noindent {\footnotesize
            {\small \em Keywords.} {\sc #1}}.}

\newcommand{\ams}[2]{  \noindent {\footnotesize
             {\small \em AMS {\rm 2000} subject classifications.
             {\rm Primary {\sc #1}; secondary {\sc #2}}}.} }

             \title*{A simple proof of Duquesne's theorem  on contour processes of conditioned Galton-Watson
trees}
\titlerunning{Contour processes of conditioned Galton-Watson trees}
             \author {Igor Kortchemski}
             \institute{Igor Kortchemski \at Universit Paris-Sud, Orsay, France \email igor.kortchemski@normalesup.org}
\begin{document}

\maketitle

\abstract{We give a simple new proof of a theorem of Duquesne, stating that
the properly rescaled contour function of a critical aperiodic
Galton-Watson tree, whose offspring distribution is in the domain of
attraction of a stable law of index $\theta \in (1,2]$, conditioned
on having total progeny $n$, converges in the functional sense to
the normalized excursion of the continuous-time height function of a
strictly stable spectrally positive Lvy process of index $\theta$.
To this end, we generalize an idea of Le Gall which consists in
using an absolute continuity relation between the conditional
probability of having total progeny exactly $n$ and the conditional
probability of having total progeny at least $n$. This new method is
robust and can be adapted to establish invariance theorems for
Galton-Watson
 trees having $n$ vertices whose degrees are prescribed to belong to a fixed subset of the
 positive integers.}
 
 \bigskip
 

 \bigskip

\ams{60J80,60F17,60G52}{05C05}

\section*{Introduction}

In this article, we are interested in the asymptotic behavior of
critical Galton-Watson trees whose offspring distribution may have
infinite variance. Aldous \cite{A2} studied the shape of large
critical Galton-Watson trees whose offspring distribution has finite
variance and proved that their properly rescaled contour functions
converge in distribution in the functional sense to the Brownian
excursion. This seminal result has motivated the study of the
convergence of other rescaled paths obtained from Galton-Watson
trees, such as the Lukasiewicz path (also known as the Harris walk)
and the height function. In \cite{MM}, under an additional
exponential moment condition, Marckert \& Mokkadem showed that the
rescaled Lukasiewicz path, height function and contour function all
converge in distribution to the same Brownian excursion. In
parallel, unconditional versions of Aldous' result have been
obtained in full generality. More precisely, when the offspring
distribution is in the domain of attraction of a stable law of index
$\theta \in(1,2]$, Duquesne \& Le Gall \cite{DuquesneLG} showed
that the concatenation of rescaled Lukasiewicz paths of a sequence
of independent Galton-Watson trees converges in distribution to a
strictly stable spectrally positive Lvy process $X$ of index
$\theta$, and the concatenation of the associated rescaled height
functions (or of the rescaled contour functions) converges in distribution to the so-called continuous-time
height function associated with $X$. In the same monograph, Duquesne \& Le Gall explained how to deduce a limit theorem for
Galton-Watson trees conditioned on having at least $n$ vertices from
the unconditional limit theorem. Finally, still in the stable
case, Duquesne \cite{Duquesne} showed that  the rescaled
Lukasiewicz path of a Galton-Watson tree conditioned on having $n$
vertices converges in distribution to the normalized excursion of
the Lvy process $X$ (thus extending Marckert \& Mokkadem's result)
and that the rescaled height and contour functions of a Galton-Watson
tree conditioned on having $n$ vertices converge in distribution to
the normalized excursion of the continuous-time height function $\H$
associated with $X$ (thus extending Aldous' result).

\bigskip

In this work, we give an alternative proof of  Duquesne's result, which is based on an
idea that appeared in the recent papers \cite{LGIto,LGM}. Let us
explain our approach after introducing some notation. For every $ x \in \mathbb{R}$, let $\fl{x}$ denote the greatest integer smaller than or equal to $x$. If $I$ is an interval, let $ \C(I, \R)$ be the space of all continuous functions $I \to \R$ equipped with the topology of uniform convergence on compact subsets of $I$. We also let $\D(I , \R)$ be the space of all right-continuous with left limits (cdlg) functions $I \to \R$,
endowed with the Skorokhod $J_1$-topology (see \cite[chap. 3]{Bill}, \cite[chap. VI]{Shir} for background concerning the Skorokhod topology). Denote by $\P_\mu$
the law of the Galton-Watson tree with offspring distribution $\mu$.
The total progeny of a tree $\tau$ will be denoted by $\zt$. Fix
$\theta \in (1,2]$ and let $(X_t)_{t \geq 0}$ be the spectrally
positive Lvy process with Laplace exponent $ \E[ \exp(- \lambda
X_t)]= \exp(t \lambda ^ \theta)$.
\begin{enumerate}
\item[(0)] We fix a critical offspring distribution $\mu$ in the domain of
attraction of a stable law of index $\theta \in (1,2]$. If $U_1,U_2, \ldots$ are i.i.d. random variables with distribution $\mu$, and $W_n=U_1+ \cdots+U_n-n$,
there exist positive constants $(B_n)_{n \geq 0}$ such that
$W_n/B_n$ converges in distribution to $X_1$.
\item[(i)] Fix $a \in (0,1)$. To simplify notation, set $ \W^{a,(n)}=( \W^{a,(n)}_j, 0 \leq j \leq \fl{na})$ where
$ \W^{a,(n)}_j= \W_{j}(\tau)/B_n$
and $ \W(\tau)$ is the Lukasiewicz path of $\tau$ (see Section \ref {sec:coding} below for its definition). Then for every function  $f_n : \Z^{\fl{an}+1} \rightarrow \R_+$,
the following absolute continuity relation holds:
\begin{eqnarray}\label{cv:intro}
\Esmu{f_n( \W^{a,(n)})|\, \zt =n}  = \Esmu{f_n( \W^{a,(n)})
D^{(a)}_n\left( \W_{\fl{an}} ( \tau)\right) |\, \zt \geq n}
\end{eqnarray}
with a certain function $D^{(a)}_n: \{ -1,0,1,2, \ldots\} \rightarrow \R_+$.
\item[(ii)] We establish the existence of a measurable function $\Gamma_a: \R_+ \rightarrow \R_+$ such
that the quantity $\left|D^{(a)}_n(j) - \Gamma_a(j/B_n)\right|$ goes to $0$ as $n \rightarrow \infty$,
uniformly in values of $j$ such that $j/B_n$ stays in a compact
subset of $\R^*_{+}$. Furthermore, if $H$ denotes the continuous-time height
process associated with $X$ and $\Nn$ stands for the It excursion measure
of $X$ above its infimum, we have for every bounded measurable function $F :
\D([0,a],\R) \rightarrow \R_+$: \begin{equation}\label{eq:c}\Nn
\left( F( (H_t)_{0 \leq t \leq a }) \Gamma_a(X_a) | \, \z>1 \right)
=\Nn \left( F( (H_t)_{0 \leq t \leq a }) | \, \z=1
\right),\end{equation} where $ \zeta$ is the duration of the excursion under $ \Nn$.
\item[(iii)] We show that under $\Pmu[\, \cdot \, | \zt=n]$, the rescaled height function converges in distribution
 on $[0,a]$ for every $a \in (0,1)$. To this end, we fix a bounded continuous function $F: \D([0,a],\R) \rightarrow \mathbb{R}_+$ and apply formula (\ref{cv:intro})
 with
$f_n( \W^{a,(n)})=F\left(\frac{B_n}{n} H_{ \fl{n t}}(\tau);0\leq t\leq
a\right)$ where $H(\tau)$ is the height function of the tree $\tau$.
Using the previously mentioned result of Duquesne \& Le Gall concerning
Galton-Watson trees conditioned on having at least $n$ vertices, we
show that we can restrict ourselves to the case where
$ \W_{\fl{an}} ( \tau)/B_n$ stays in a compact subset of $\R^*_+$, so that we
can apply (ii) and obtain that:
\begin{eqnarray*} &&\lim_{n \rightarrow \infty}
\Esmu{ \left.  F\left(\frac{B_n}{n} H_{ \fl{n t}}(\tau);\, 0 \leq t \leq a \right)
\right | \zt=n} \\
&& \qquad\qquad\qquad\qquad = \lim_ {n \rightarrow \infty} \Esmuzgeqn{ F\left(\frac{B_n}{n} H_{ \fl{n t}}(\tau);\, 0 \leq t \leq a \right)D^{(a)}_n\left( \W_{\fl{an}} ( \tau) \right)}\\
&& \qquad\qquad\qquad\qquad  = \Nn \left( F( H_t; 0 \leq t \leq a) \Gamma_a(X_a) \, \|
\, \z
> 1\right)\\
&& \qquad\qquad\qquad\qquad =  \Nn \left( F( H_t; 0 \leq t \leq a )\, \| \, \z =
1\right).\end{eqnarray*}
\item[(iv)] By using a relationship between the contour function and
the height function which was noticed by Duquesne \& Le Gall in
\cite{DuquesneLG}, we get that, under $\Pmu[\, \cdot \, |
\zt=n]$, the scaled contour function converges in distribution
 on $[0,a]$.
 \item[(v)] By using the time reversal invariance property of the contour function, we
 deduce that under $\Pmu[\, \cdot \, | \zt=n]$,
the scaled contour function converges in distribution
 on the whole segment $[0,1]$.
 \item[(vi)] Using once again the relationship between the contour function and
the height function, we deduce that, under $\Pmu[\, \cdot \, |
\zt=n]$, the scaled height function converges in distribution
 on $[0,1]$.
\end{enumerate}

In the case where the variance of $\mu$
is finite, Le Gall gave an alternative proof of Aldous' theorem in
\cite[Theorem 6.1]{LGIto} using a similar approach based on a strong local limit theorem. There are additional difficulties in the infinite variance case, since no such
theorem is known in this case. 

 Let us finally discuss the advantage of this new method.
 Firstly, the proof is simpler and less technical.
 Secondly, we believe that this approach is robust and can be adapted to other situations.
 For instance, using the same ideas, we have established invariance theorems for Galton-Watson
 trees having $n$ vertices whose degrees are prescribed to belong to a fixed subset of the
 nonnegative integers \cite{Kleaves}.

\bigskip

The rest of this text is organized as follows. In Section 1, we
present the discrete framework by defining Galton-Watson trees and
their codings. We explain how the local limit theorem gives
information on the asymptotic behavior of large Galton-Watson trees
and present the discrete absolute continuity relation appearing in
(\ref{cv:intro}). In Section 2, we discuss the continuous framework:
we introduce the strictly stable spectrally positive Lvy process,
its It excursion measure $\Nn$  and the associated
continuous-time height process. We also prove the absolute
continuity relation (\ref{eq:1}). Finally, in Section 3 we give
the new proof of Duquesne's theorem by carrying out steps (i-vi).

\bigskip

\textbf{Acknowledgments.} I am deeply indebted to Jean-Franois Le Gall for insightful
discussions and for making many useful suggestions on first versions of this manuscript.
\bigskip

\textbf{Notation and main assumptions.} Throughout this work $\theta \in (1,2]$ is a fixed parameter. We
consider a probability distribution $(\mu(j))_{j \geq 0}$ on the
nonnegative integers satisfying the following three conditions:
\begin{enumerate}
\item[(i)] $\mu$ is critical, meaning that $\sum_{k=0}^\infty k \mu(k)=1$.
\item[(ii)] $\mu
$ is in the domain of attraction of a stable law of index $\theta
\in (1,2]$. This
means that either the variance of $\mu$ is positive and finite, or $\mu([j,\infty))= j^{-\theta} L(j)$,
where $L: \R_+ \rightarrow \R_+$ is a function such that $L(x)>0$ for $x$ large enough and $\lim_{x
\rightarrow \infty} L(tx)/L(x)=1$ for all $t>0$ (such a function is
called slowly varying). We refer to \cite{Bingham} or \cite[chapter
3.7]{Durrett} for details.
\item[(iii)] $\mu$ is aperiodic, which means that
the additive subgroup of the integers $\mathbb{Z}$ spanned by $\{j;
\, \mu(j) \neq 0 \}$ is not a proper subgroup of $\mathbb{Z}$.
\end{enumerate}
We introduce condition (iii) to avoid unnecessary complications, but
our results can be extended to the periodic case.

In what follows, $(X_t)_{t \geq 0}$ will stand for the spectrally positive Lvy process with Laplace exponent $ \E[ \exp(- \lambda X_t)]=\exp(t \lambda ^ \theta)$ where $t, \lambda \geq 0$ and $p_1$ will denote the density of $X_1$.
Finally, $\nu$ will stand for the probability measure on $ \Z$
defined by $\nu(k)=\mu(k+1)$ for $k \geq -1$. Note that $\nu$ has
zero mean.

\section{The discrete setting : Galton-Watson trees}
\subsection{Galton-Watson trees}

\begin{definition}Let $\N=\{0,1,\ldots\}$ be the set of all nonnegative integers and $\N^*=\{1,\ldots\}$.  Let also $U$ be the set of all labels:
$$U=\bigcup_{n=0}^{\infty} (\N^*)^n,$$
where by convention $(\N^*)^0=\{\emptyset\}$. An element of $U$ is a
sequence $u=u_1 \cdots u_j$ of positive integers, and we set
$|u|=j$, which represents the \og generation \fg \, of $u$. If $u=u_1
\cdots u_j$ and $v=v_1 \cdots v_k$ belong to $U$, we write $uv=u_1
\cdots u_j v_1 \cdots v_k$ for the concatenation of $u$ and $v$. In
particular, note that $u \emptyset=\emptyset u = u$. Finally, a
\emph{rooted ordered tree} $\tau$ is a finite subset of $U$ such
that:
\begin{itemize}
\item[1.] $\emptyset \in \tau$,
\item[2.] if $v \in \tau$ and $v=uj$ for some $j \in \N^*$, then $u
\in \tau$,
\item[3.] for every $u \in \tau$, there exists an integer $k_u(\tau)
\geq 0$ such that, for every $j \in \N^*$, $uj \in \tau$ if and only
if $1 \leq j \leq k_u(\tau)$.
\end{itemize}
In the following, by \emph{tree} we will mean rooted ordered
tree. The set of all trees is denoted by $\T$. We will often view each vertex of
a tree $\tau$ as an individual of a population whose $\tau$ is the
genealogical tree. The total progeny of $\tau$ will be
denoted by $\zeta(\tau)=\Card(\tau)$. Finally, if $\tau$ is a tree and $u \in \tau$, we
set $T_u \tau=\{v \in U; \, uv \in
\tau\}$, which is itself a tree.
\end{definition}

\begin{definition}Let $\rho$ be a probability measure on $\N$ with mean less than or equal to $1$ and such that $\rho(1)<1$. The law of the
Galton-Watson tree with offspring distribution $\rho$ is the unique
probability measure $\P_\rho$ on $\T$ such that:
\begin{itemize}
\item[1.] $\P_\rho[k_\emptyset=j]=\rho(j)$ for $j \geq 0$,
\item[2.] for every $j \geq 1$ with $\rho(j)>0$, conditionally on $\{k_ \emptyset=j\}$, the shifted trees
$T_1 \tau, \ldots, T_j \tau$ are i.i.d. with distribution $\P_\rho$.
\end{itemize}
A random tree whose distribution is $\P_\rho$ will be called a
$\GW_\rho$ tree.
\end{definition}

\subsection{Coding Galton-Watson trees}
\label {sec:coding}

We now explain how trees can be coded by three different functions.
These codings are crucial in the understanding of large
Galton-Watson trees.

\begin{figure*}[h!]
 \begin{minipage}[c]{0.45 \linewidth}
   \centering
      \includegraphics[scale=0.5]{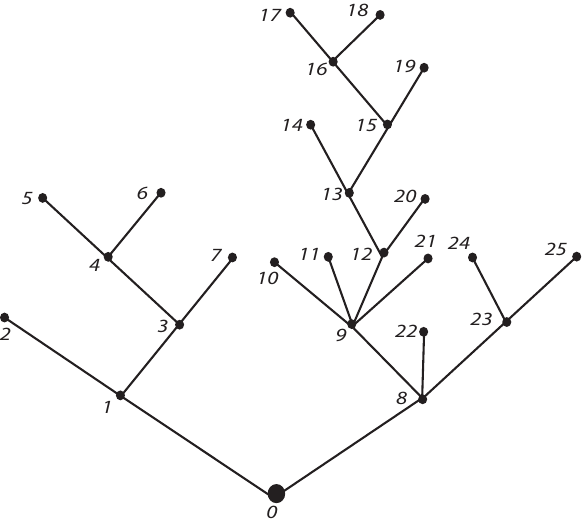}
   \end{minipage}
   \begin{minipage}[c]{0.45 \linewidth}
   \centering
      \includegraphics[scale=0.4]{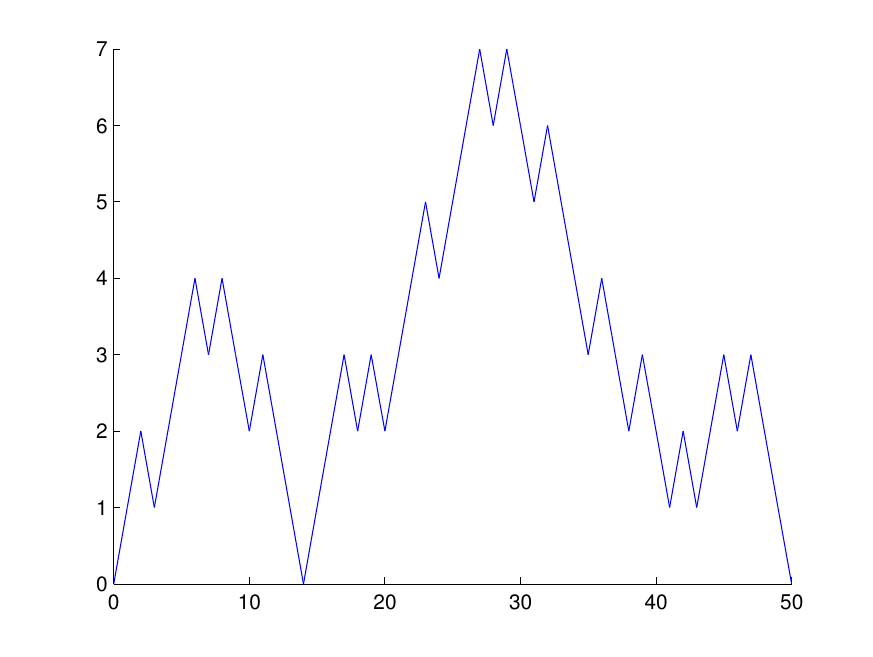}
   \end{minipage}
   \caption{\label{fig:tree1}A tree $\tau$ with its vertices indexed in
lexicographical order and its  contour function $(C_{u}(\tau);\, 0
\leq u \leq 2(\zt-1)$. Here, $\zeta(\tau)=26$.}
   \begin{minipage}[c]{0.45 \linewidth}
   \centering
      \includegraphics[scale=0.4]{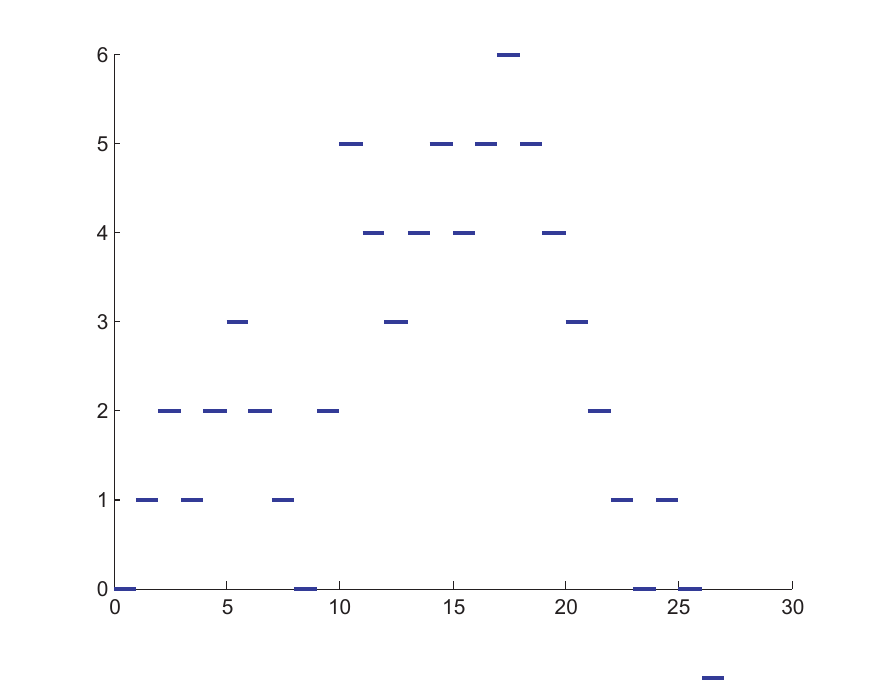}
   \end{minipage}
   \begin{minipage}[c]{0.45 \linewidth}
   \centering
      \includegraphics[scale=0.4]{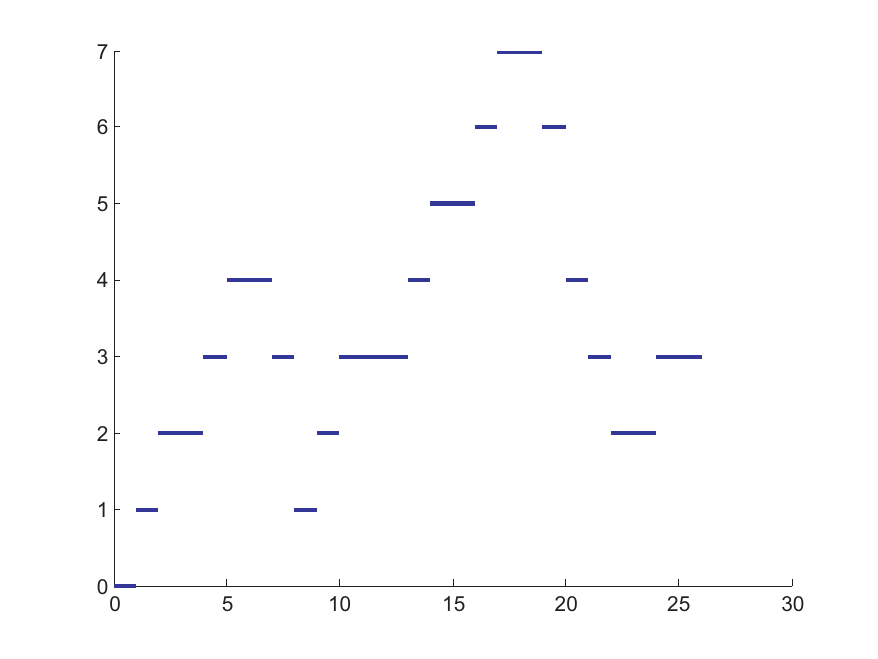}
   \end{minipage}
   \caption{\label{fig:tree} The Lukasiewicz path $( \W_{ \fl{u}}(\tau);\, 0 \leq u < \zt+1)$ and the height function
$(H_{u}(\tau), 0 \leq u \leq \zt)$ of the tree $\tau$ of Fig. \ref{fig:tree1}.}
\end{figure*}

\begin{definition}\label{def:fonctions}
We write $u<v$ for the lexicographical order on the labels $U$ (for
example $\emptyset<1<21<22$). Consider a tree $\tau$ and order the
individuals of $\tau$ in lexicographical order:
$\emptyset=u(0)<u(1)<\cdots<u(\zeta(\tau)-1)$. The height process
$H(\tau)=(H_n(\tau), 0 \leq n < \zeta(\tau))$ is defined, for $0
\leq n < \zeta(\tau)$, by $H_n(\tau)=|u(n)|$. For technical
reasons, we set $H_{k}(\tau)=0$ for $k \geq \z(\tau)$. We extend $H(\tau)$ to $ \R_+$ by linear interpolation
by setting $H_t (\tau) = (1 - \{t\})H_{ \fl{t}}(\tau) + \{t\}H_{\fl{t}+1}(\tau)$ for $0 \leq t \leq \zeta( \tau )$,
where $\{t\} = t-\fl{t}$.

 Consider a particle that starts from
the root and visits continuously all edges at unit speed (assuming
that every edge has unit length), going backwards as little as
possible and respecting the lexicographical order of vertices. For
$0 \leq t \leq 2(\zt-1)$, $C_t(\tau)$ is defined as the distance to
the root of the position of the particle at time $t$. For technical
reasons, we set $C_t(\tau)=0$ for $t \in [2(\zt-1), 2 \zt]$. The
function $C(\tau)$ is called the contour function of the tree
$\tau$. See Figure \ref{fig:tree1} for an example, and \cite[Section
2]{Duquesne} for a rigorous definition.

Finally, the Lukasiewicz path $ \W(\tau)=( \W_n(\tau), n \geq 0 )$ of a tree $\tau$ is defined by $ \W_0(\tau)=0$, $ \W_{n+1}(\tau)= \W_{n}(\tau)+k_{u(n)}(\tau)-1$
for $0 \leq n \leq \zeta(\tau)-1$ and $ \W_k( \tau)=0$ for $ k> \zt$. For $u \geq 0$, we set $ \W_u(\tau)=  \W_ {\fl{u}}(\tau)$.
\end{definition}
Note that necessarily $ \W_{\zt}(\tau)=-1$.

\bigskip

Let $(W_n; n \geq 0)$ be a
random walk which starts at $0$ with jump distribution $\nu(k)=\mu(k+1)$ for $k \geq -1$. For $j \geq 1$, define $\zeta_j =\inf\{n \geq 0; \, W_n=-j\}$.

\begin{proposition}\label{prop:RW} $(W_0,W_1,\ldots,W_{\zeta_1})$ has the same distribution as the Lukasiewicz path of a
$\GW_\mu$ tree. In particular, the total progeny of a
$\GW_\mu$ tree has the same law as $\z_1$.\end{proposition}

\begin{proof}See \cite[Proposition 1.5]{RandomTrees}.\end{proof}

We will also use the following well-known fact (see e.g. Lemma 6.1 in \cite{Pitman} and the discussion that follows). 

\begin{proposition}\label{prop:kemperman} For every integers $ 1 \leq j \leq n$, we have $\P[ \z_j =n] = \frac{j}{n} \P[W_n=-j]$.\end{proposition}

\subsection{Slowly varying functions}

Slowly varying functions appear in the study of domains of
attractions of stable laws. Here we recall some properties of these
functions in view of future use.

Recall that a positive measurable function $L: \R_+ \rightarrow
\R_+$ is said to be slowly varying if $ L(x)>0$ for $x$ large enough and, for all $t>0$, $ L(tx)/L(x)
\rightarrow 1$ as $x \rightarrow \infty$. A useful result concerning
these functions is the so-called Representation Theorem, which
states that a function $L: \R_+ \rightarrow \R_+$ is slowly varying
if and only if it can be written in the form:
$$L(x)=c(x) \exp \left( \int_1^x \frac{\e(u)}{u} du\right), \qquad x \geq 0,$$
where $c$ is a nonnegative measurable function having a finite positive
limit at infinity and $\e$ is a measurable function tending to $0$
at infinity. See e.g. \cite[Theorem 1.3.1]{Bingham} for a proof. The
following result is then an easy consequence.

\begin{proposition}\label{prop:slow}Fix $\epsilon>0$ and let  $L: \R_+ \rightarrow \R$ be a
slowly varying function. There exist two constants $C>1$ and $N>0$ such that
$\frac{1}{C}x^{-\e}\leq{L(nx)}/{L(n)} \leq Cx^\e$ for every integer
$n\geq N$ and $x \geq 1$.
\end{proposition}

\subsection{The Local Limit Theorem}

\begin{definition}A subset $A \subset \Z$ is said to be lattice if there exist $b \in \Z$ and an integer $d  \geq 2$ such that
$A \subset  b+ d \Z$. The largest $d$ for which this statement holds
is called the span of $A$. A measure on $\Z$ is said to be lattice
if its support is lattice, and a random variable is said to be
lattice if its law is lattice.\end{definition}

\begin{remark}\label{rem:non-lattice}Since $\mu$ is supposed to be critical and aperiodic, using the fact that
$\mu(0)>0$, it is an exercise to check that the probability measure
$\nu$ is non-lattice.\end{remark}

Recall that $p_1$ is the density of $X_1$. It is well known that $p_1(0)>0$, that $p_1$ is positive, bounded and continuous, and that the absolute value of the
derivative of $p_1$ is bounded over $\R$ (see e.g. \cite[I. 4.]{Zolotarev}). The following theorem will allow us to find estimates for the
probabilities appearing in Proposition \ref{prop:kemperman}.

\begin{theorem}[Local Limit Theorem]\label{thm:locallimit}
Let $(Y_n)_{n \geq 0}$ be a random walk on $\Z$ started from $0$
such that its jump distribution is in the domain of attraction of a
stable law of index $\theta \in (1,2]$. Assume that $Y_1$ is
non-lattice, that $\E[Y_1]=0$ and that $Y_1$ takes values in $\N \cup \{-1\}$.
\begin{enumerate}
\item[(i)] There exists an increasing sequence of positive real numbers $(a_n)_{n
\geq 1}$ such that $Y_n/a_n$ converges in distribution to $X_1$.
\item[(ii)] We have $ \d \lim_{n \rightarrow \infty} \sup_{k \in \Z}
\left| a_n \P[Y_n=k]-p_1\left( \frac{k}{a_n}\right)
\right|=0.$
\item[(iii)] There exists a
slowly varying function $l: \R_+ \rightarrow \R_+$ such that
$a_n=n^{1/\theta} l(n)$.
\end{enumerate}
\end{theorem}

\begin{proof} For (i), see \cite[Section XVII.5, Theorem 3]{Feller}
and \cite[Section 8.4]{Bingham}. The fact that $(a_n)$ may be chosen
to be increasing follows from  \cite[Formula 3.7.2]{Durrett}. For
(ii), see \cite[Theorem 4.2.1]{IL}. For (iii), it is shown in
\cite[p. 46]{IL} that ${a_{kn}}/{a_n}$ converges to $k^{1/\theta}$
for every integer $k \geq 1$. Since $(a_n)$ is increasing, a theorem of de Haan (see \cite[Theorem 1.10.7]{Bingham}) implies that there exists a
slowly varying function $l: \R_+ \rightarrow \R_+$ such that
$a_n=l(n) n^{1/\theta}$ for every positive integer $n$.
\end{proof}

Let $(W_n)_{n \geq 0}$ be as in Proposition \ref {prop:RW} a random walk started from $0$ with jump
distribution $\nu$. Since $\mu$ is in the domain of attraction of a
stable law of index $\theta$, it follows that $\nu$ is in the same
domain of attraction, and $W_1$ is not lattice by Remark
\ref{rem:non-lattice}. Since $\nu$ has zero mean, by the preceding
theorem there exists an increasing sequence of positive integers
$(B_n)_{n \geq 1}$ such that $B_n \rightarrow \infty$ and $W_n/B_n$ converges in distribution towards $X_1$ as $n \rightarrow \infty$.
In what follows, the sequence $(B_n)_{n \geq 1}$ will be fixed, and $h: \R_+ \rightarrow \R_+$
will stand for a slowly varying function such that $B_n=h(n)n
^{1/\theta}$.

\begin{lemma}\label{lem:locallimit}We have:
$$ (i) \quad\Prmu{\zeta(\tau)=n} \quad \mathop{\sim}_{n \rightarrow \infty} \quad \frac{ p_1(0) }{n^{1/\theta+1} h(n)}, \qquad \qquad (ii) \quad\Prmu{\zeta(\tau) \geq n} \quad \mathop{\sim}_{n \rightarrow \infty} \quad \frac{\theta p_1(0)}{n^{1/\theta} h(n)},$$
where we write $a_n \sim b_n$ if $a_n/b_n \rightarrow 1$.
\end{lemma}

\begin{proof}We keep the notation of Proposition \ref{prop:RW}. Proposition \ref{prop:kemperman} gives that:
\begin{equation}\label{eq:1}\P_\mu[\zeta(\tau)=n]=\frac{1}{n}
\P[W_n=-1].\end{equation} For (i), it suffices to notice that the
local limit theorem (Theorem \ref{thm:locallimit}) and the continuity of $p_1$ entail
$\P[W_n=-1] \sim {p_1(0)}/({h(n)n^{1/\theta}})$.
For (ii), we use (i) to write:
$$
\Prmu{\zeta(\tau) \geq n} = \sum_{k=n}^{\infty}
\left( \frac{1}{h(k) k^{1+1/\theta} }
p_1\left(0 \right) + \frac{1}{h(k)
k^{1+1/\theta} } \delta(k)\right),
$$
where $\delta(k)$ tends to $0$ as $k \rightarrow \infty$. We can rewrite this in the form:
\begin{equation}
\label{eq:somintegrale}h(n) n^{1/\theta} \Prmu{\zeta(\tau) \geq n}  =  \int_1 ^\infty du f_n(u), \end{equation}
where:
$$f_n(u)=  \frac{h(n) n^{1/\theta+1}}{h(\fl{n u}) \fl{n u}^{1+1/\theta} }
 \left( p_1\left(0\right)+ \delta(\fl{nu}) \right).$$
 For fixed $u \geq 1$,  $f_n(u)$ tends to $\frac{p_1(0)}{
 u^{1/\theta+1}}$ as $n \rightarrow \infty$. Choose $\e \in (0,1/\theta)$. By
 Proposition
 \ref{prop:slow}, for every sufficiently large positive integer $n$  we have $f_n(u) \leq C/u^{1+1/\theta-\e}$ for every $u \geq 1$, where $C$ is a positive constant.
 The dominated convergence theorem allows us to infer that:
 $$\lim_{n \rightarrow \infty} \int_1 ^\infty du f_n(u) = \int_{1}^{\infty}du \frac{p_1(0)}{
 u^{1/\theta+1}}= \theta p_1(0),$$
 and the desired result follows from \eqref{eq:somintegrale}.
\end{proof}

\subsection{Discrete absolute continuity}

The next lemma in another important ingredient of our approach.

\begin{lemma}[Le Gall \& Miermont]\label{lem:AC}
Fix $a \in (0,1)$. Then, with the notation of Proposition \ref{prop:kemperman}, for every $n \geq 1$ and for every bounded nonnegative function $f_n$ on $\Z^{\fl{an}+1}$:
\begin{equation}
\label{eq:rel2}
\Es{f_n(W_0, \ldots, W_ {\fl{na}})|\, \z_1 =n}=
\Es{ \left. f_n(W_0, \ldots, W_ {\fl{na}} ) \frac{\phi_{n-[a
n]}(W_{\fl{an}}+1)/\phi_n(1)}{\phi^*_{n-\fl{an}}(W_{ \fl{a n}}+1)/\phi^*_n(1)} \right|\, \z_1 \geq n}, 
\end{equation}
where $ \phi_p(j)=\Pr{\z_j=p}$ and $
\phi^*_p(j)=\Pr{\z_j \geq p}$ for every integers $j \geq 1$ and $ p \geq 1$.
\end{lemma}

\begin{proof}
This result follows from an application of the Markov property to the random walk $W$ at time  $\fl{an}$. See \cite[Lemma 10]{LGM} for details in a slightly different setting.
\end{proof}

\section{The continuous setting: stable Lvy processes}

\subsection{The normalized excursion of the Lvy process and the continuous-time height process}
\label {subsec:ito}

We follow the presentation of \cite{Duquesne}. The underlying
probability space will be denoted by $(\Omega, \mathcal{F}, \P)$. Recall that $X$ is a strictly stable
spectrally positive Lvy process with index $\theta \in (1,2]$ such
that for $\lambda>0$:
\begin{equation}\label{eq:X}\E[\exp(-\lambda X_t)]=\exp(t
\lambda^\theta).\end{equation} We denote
the canonical filtration generated by $X$ and augmented with the
$\P$-negligible sets  by $(\F_t)_{t \geq 0}$. See \cite{Bertoin} for the proofs of
the general assertions of this subsection concerning Lvy processes.
In particular, for $\theta=2$ the process $X$ is $\sqrt{2}$ times
the standard Brownian motion on the line. Recall that $X$ has the
following scaling property: for $c>0$, the process $(c^{-1/\theta}
X_{ct}, t \geq 0$) has the same law as $X$. In particular, the density $p_t$ of the law of $X_t$ enjoys the following scaling property:
\begin{equation}\label{eq:scaling}p_t(x)=t^{-1/\theta} p_1(x t^{-1/\theta})\end{equation}
for $x \in \R$, $t > 0$. The following notation will be useful: for $s<t$, we set $I^s_t = \inf_{[s,t]} X$ and $I_t=\i{0,t} X$.
Notice that the process $I$ is continuous since $X$ has no negative
jumps.

The process $X-I$ is a strong Markov process and $0$ is regular for
itself with respect to $X-I$. We may and will choose $-I$ as the
local time of $X-I$ at level $0$. Let $(g_i,d_i), i \in \mathcal{I}$
be the excursion intervals of $X-I$ above $0$. For every $i \in \mathcal{I}$ and $ s \geq  0$, set $\omega_s^i= X_{(g_i+s) \wedge d_i}-X_{g_i}$. We view $ \omega^i$ as an element of the excursion space $ \mE$, which is defined by:
$$ \mE=  \{ \omega \in \D(\R_+, \R); \, \omega(0)=0 \textrm{ and } \z( \omega):=
 \sup\{s>0; \omega(s)>0\,\} \in (0, \infty)  \}.$$

From It's excursion theory, the point measure
$$\mathcal{N}( dt d \omega)= \sum_{i \in \mathcal{I}} \delta_{(-I_{g_i},\omega^i)}$$
is a Poisson measure on $ \R_+ \times \mE$ with intensity $dt \Nn(d\omega)$, where $\Nn(d
\omega)$ is a $\sigma$-finite measure on $\mE$. By
classical results, $ \Nn(\zeta>t)={\Gamma(1-1/\theta)}^ {-1}t^{-1/\theta}$.
 Without risk of confusion, we will also use the
notation $X$ for the canonical process on the space $\D(\R_+, \R)$.

We now define the normalized excursion of $X$. Let us first recall the It description of the excursion measure (see \cite{Chaumont}
or \cite[Chapter VIII.4]{Bertoin} for details). Define for $\lambda>0$ the re-scaling operator $S^{(\lambda)}$
on $ \mE$ by $S^{(\lambda)}(\omega)=\left( \lambda^{1/\theta} \omega(s/\lambda),s \geq 0 \right).$
Then there exists a unique collection of probability measures
$(\Nn_{(a)}, a>0)$ on $\mE$ such that the
following properties hold.
\begin{enumerate}
\item[(i)] For every $a>0$, $\Nn_{(a)}(\z=a)=1$.
\item[(ii)] For every $\l>0$ and $a>0$, we have
$S^{(\l)}(\Nn_{(a)})=\Nn_{(\l a)}$.
\item[(iii)] For every measurable subset $A$ of $ \mE$: $ \d \Nn(A)= \int_0 ^\infty \Nn_{(a)}(A) \frac{da}{\theta \Gamma(1-1/\theta)
a^{1/\theta+1}}$.
\end{enumerate}
The probability distribution $\Nn_{(1)}$ on cdlg paths with unit
lifetime is called the law of the normalized excursion of $X$ and will sometimes be
denoted by $\Nn( \, \cdot \, | \z=1)$.  In particular, for
$\theta=2$, $\Nn_{(1)}$ is the law of $ \sqrt {2}$ times the normalized
excursion of linear Brownian motion. Informally,  $\Nn( \, \cdot \,
| \z=1)$ is the law of an excursion
conditioned to have unit lifetime.

We will also use the so-called continuous-time height process $H$
associated with $X$ which was introduced in \cite{LeJan}. If
$\theta=2$, $H$ is set to be equal to $X-I$. If $\theta \in
(1,2)$, the process $H$ is defined for every $t \geq 0$ by:
$$H_t:=\lim_{\e \rightarrow 0} \frac{1}{\e} \int_0 ^t
\mathbbm{1}_{\{X_s<I^s_t+\e\}}ds,$$ where the limit exists in $\P$-probability and in
$\Nn$-measure on $  \{t < \zeta\}$. The definition of $H$ thus makes sense under $ \P$ or under $\Nn$. The process $H$ has a
continuous modification both under $ \P$ and under $ \Nn$ (see \cite[Chapter 1]{DuquesneLG} for
details), and from now on we consider only this modification. Using simple scale arguments one can also define $H$ as a continuous random process under $\Nn( \, \cdot \, | \z=1)$.  For our purposes, we will need the fact that, for every
$a \geq 0$, $(H_t)_{0 \leq t \leq a}$ is a measurable function of
$(X_t)_{0 \leq t \leq a}$.

\subsection{Absolute continuity property of the It measure}

We now present the continuous counterpart of the discrete absolute
continuity property appearing in Lemma \ref{lem:AC}. We follow the
presentation of \cite{LGIto} but generalize it to the stable case. The following proposition is classical (see e.g. the proof of Theorem 4.1 in \cite[Chapter XII]{RevuzYor},
which establishes the result for Brownian motion).

\begin{proposition}\label{prop:markovian}Fix $t>0$. Under the conditional probability measure $\Nn( \, \cdot \, | \z>t)$, the process
$(X_{t+s})_{s \geq 0}$ is Markovian with the transition kernels of a
strictly stable spectrally positive Lvy process of index $\theta$
stopped upon hitting $0$.\end{proposition}

We will also use the following result (see \cite[Corollary 2.3]{BertoinSub} for a proof).
\begin{proposition}\label{prop:passagetime}Set $q_s(x)=\frac{x}{s}p_s(-x)$ for $x,s>0$. For $x\geq 0$, let
$T(x)= \inf \{ t \geq 0; \, X_t<-x\}$ be the first passage time of
$-X$ above $x$. Then $\Pr{T(x) \in dt}= q_t(x) dt$ for every $x >
0$.\end{proposition}

Note that $q_s$ is a positive continuous function on $(0, \infty)$, for every $s>0$. It is also known that $q_s$ is bounded by a constant which is uniform when $s$ varies over $[ \e, \infty)$, $ \epsilon>0$ (this follows from e.g. \cite[I. 4.]{Zolotarev}).

\begin{proposition}\label{prop:Gamma}For every $a \in (0,1)$ and $x>0$ define:
$$\Gamma_a(x)=\frac{\theta q_{1-a}(x)}{\int_{1-a} ^ \infty ds \, q_s(x)}.$$
Then for every measurable bounded function $G : \D([0,a],\R^2) \rightarrow \R_+$:
$$\Nn \left( G( (X_t)_{0 \leq t \leq a },(H_t)_{0 \leq t \leq a }) \Gamma_a(X_a) | \, \z>1 \right)
=\Nn \left( G( (X_t)_{0 \leq t \leq a }, (H_t)_{0 \leq t \leq a }) | \, \z=1 \right).$$\end{proposition}

\begin{proof}Since $(H_t)_{0 \leq t \leq a }$ is a measurable function of
$(X_t)_{0 \leq t \leq a }$, it is sufficient to prove that for every bounded
measurable function $F : \D([0,a],\R) \rightarrow \R_+$:
$$\Nn \left( F( (X_t)_{0 \leq t \leq a }) \Gamma_a(X_a) | \, \z>1 \right)
=\Nn \left( F( (X_t)_{0 \leq t \leq a }) | \, \z=1 \right).$$  To this end, fix $r \in [0,a]$, let
$f,g : \R_+ \rightarrow \R_+$ be two bounded continuous
functions and let $h:  \R^*_+ \rightarrow \R_+$ be a continuous function.  Using the notation of Proposition \ref{prop:passagetime}, we have:
\begin{eqnarray}
\Nn\left(f(X_r) h(X_a) g(\z) \mathbbm{1}_{\{\z > a \}} \right)&=& \Nn\left(f(X_r)\mathbbm{1}_{\{\z > a \}}  \Es{h(x) g(a+T(x))}_ {x=X_a} \right) \notag \\
 &=&  \int_0 ^\infty ds \,
g(a+s) \Nn \left( f(X_r) h(X_a) q_s(X_a) \mathbbm{1}_{\{\z > a
\}} \right) \notag \\
&=& \int_a^ \infty du g(u) \Nn \left( f(X_r) h(X_a)q_ {u-a}(X_a) \mathbbm{1}_ {  \{ \z>a\}} \right)
, \label {eq:markov}
\end{eqnarray}
where we have used Proposition \ref{prop:markovian} in the first equality and Proposition \ref{prop:passagetime} in the second equality.
 Moreover, by property (iii) in subsection \ref {subsec:ito}:
\begin{equation}
\label{eq:itodescr} \Nn\left(f(X_r) g(\z)\mathbbm{1}_{\{\z > a
\}}\right)= \int_a ^\infty du \,\frac{g(u)}{\theta \Gamma(1-1/\theta) u^{1/\theta+1}} \cdot  \Nn_{(u)} (f(X_r)).
\end{equation}
Now observe that \eqref{eq:markov} (with $h=1$) and \eqref{eq:itodescr} hold for any bounded continuous function $g$. Since both functions $u \mapsto \Nn \left( f(X_r) q_ {u-a}(X_a) \mathbbm{1}_ {  \{ \z>a\}} \right)$ and $ u \mapsto \Nn_ {(u)} \left(f(X_r) \right)$ are easily seen to be continuous over $(a, \infty)$, it follows that for every $u>a$:
$$
\Nn\left(f(X_r)q_{u-a}(X_a) \mathbbm{1}_{\{\z > a \}}
\right)=
\frac{1}{\theta\Gamma(1-1/\theta) u^{1/\theta+1} }\Nn_{(u)}\left( f(X_r) \right).$$
In particular, for $u=1$ we get:
\begin{equation}\label{eq:N1}
\Nn\left(f(X_r)q_{1-a}(X_a) \mathbbm{1}_{\{\z > a \}}
\right)=
\frac{1}{\theta\Gamma(1-1/\theta)}\Nn_{(1)}\left( f(X_r) \right).\end{equation}

On the other hand, applying \eqref{eq:markov} with $g(x) = \mathbbm {1}_{  \{x>1\}}$ and noting that $\Nn(\z
> 1)=\frac{1}{\Gamma(1-1/\theta)}$, we get:
\begin{equation}\label{eq:N2}\Nn\left( f(X_r) h(X_a) \, | \z >1 \right)= \Gamma(1-1/\theta)\Nn\left(f(X_r) h(X_a) \mathbbm{1}_{\{\z > a
\}}\int_{1-a}^\infty ds \, q_s(X_a)\right).\end{equation} By combining (\ref{eq:N2}) and (\ref{eq:N1}) we conclude that:
$$ \Nn\left( \left.f(X_r) \frac{\theta q_{1-a}(X_a)}{\int_{1-a} ^ \infty ds \, q_s(X_a)}
\, \right|\, \z >1 \right)=\Nn_{(1)}( f(X_r)).$$ One similarly shows that for $0 \leq r_1 <
\cdots < r_n \leq a$ and $f_1,\ldots,f_n: \R_+ \rightarrow \R_+$
 continuous bounded functions:
$$ \Nn\left( \left. f_1(X_{r_1}) \cdots f_n(X_{r_n}) \frac{\theta q_{1-a}(X_a)}{\int_{1-a} ^ \infty ds \, q_s(X_a)}
\, \right|\, \z >1 \right)=\Nn_{(1)}( f_1(X_{r_1}) \cdots f_n(X_{r_n})).$$  The desired result follows since the Borel
$\sigma$-field of $\D([0,a],\R)$ is generated by the coordinate
functions $X \mapsto X_r$ for $0 \leq r \leq a$ (see e.g. \cite[Theorem 12.5 (iii)]{Bill}).
\end{proof}

\section{Convergence to the stable tree}

\subsection{An invariance theorem}
We rely on the following theorem, which is similar in spirit to
Donsker's invariance theorem (see the concluding remark of
\cite[Section 2.6]{DuquesneLG} for a proof).

\begin{theorem}[Duquesne \& Le Gall]\label{thm:cvW}Let $\t_n$ be a random tree distributed according to $\Pmuzgeqn$. We have:
$$ \left( \frac{1}{B_n} W_{ \lfloor nt \rfloor}(\t_n), \frac{B_n}{n} H_{nt}(\t_n)\right)_ {t \geq 0} \quad
\mathop{\longrightarrow}^{(d)}_{n \rightarrow \infty}\quad (
X_t,H_t)_{0 \leq t \leq 1} \textrm{ under }\Nn( \, \cdot \, | \,
\z>1). $$
\end{theorem}

\subsection{Convergence of the scaled contour and height functions}

Recall the notation $\phi_n(j)=\Pr{\z_j=n}$ and $
\phi^*_n(j)=\Pr{\z_j \geq n}$.

\begin{lemma}\label{lem:technical2}
Fix $\a>0$. We have:
$$(i) \, \lim_{n \rightarrow \infty }   \sup_{1 \leq k \leq \a B_n} \left| n \phi_n(k)-
q_1\left(\frac{k}{B_{n}}\right) \right|=0, \qquad (ii) \, \lim_{n \rightarrow \infty }   \sup_{1 \leq k \leq \a B_n} \left|
\phi^*_n(k)- \int_1 ^\infty ds \, q_s\left(\frac{k}{
 B_{n} }\right) \right|=0.$$\end{lemma}

This has been proved by Le Gall in \cite{LGIto} when $\mu$ has
finite variance. In full generality, the proof is technical and is
postponed to Section 3.3.

\begin{lemma}\label{lem:cvH0a}Fix $a \in (0,1)$. Let $\t_n$ be a random tree distributed according to  $\Pmuzn$. Then the following convergence holds in distribution in the space $ \C([0,a], \R)$:
 $$\left(\frac{B_n}{n} H_{nt}(\t_n); \, 0 \leq t \leq a\right) \qquad
 \mathop{\longrightarrow}_{n \rightarrow \infty}^{(d)} \qquad
 (H_t; \, 0 \leq t \leq a)  \textrm{ under } \Nn( \, \cdot\, | \z
 =1).$$\end{lemma}

\begin{proof}Recall the notation $ \Gamma_a$ introduced in Proposition \ref {prop:Gamma}. We start by verifying that, for $\a > 1$, we have:
\begin{equation}\label{eq:approx2}\lim_{n \rightarrow \infty }
\left( \sup_{ \frac{1}{\a} B_n \leq k \leq \a B_n} \left|
\frac{\phi_{n- \fl{a n}}(k+1)/\phi_n(1)}{\phi^*_{n-\fl{a
n}}(k+1)/\phi^*_n(1)}-\Gamma_a\left(\frac{k}{B_
n}\right)\right|\right)=0.\end{equation}
To this end, we will use the existence of a constant $ \delta>0$ such that, for $n$
sufficiently large,
\begin{equation}
\label{eq:delta} \inf_ {\frac{1}{\a} B_n \leq k \leq \a B_n} \int_ {1} ^ { \infty}ds q_s \left( \frac{k+1}{B_ { n- \fl{an}}} \right)> \delta.
\end{equation}
The existence of such $ \delta$ follows from the fact that, for every $ \beta>1$, $ \inf_{\frac{1}{\b} \leq x \leq \b }   \int_ {1} ^ { \infty}ds q_s \left(x \right)>0.$
We will also need the fact that  for every  $ \beta>1$ there exists a constant $C>0$ such :
\begin{equation}
\label{eq:C}\sup_{ \frac{1}{ \beta} \leq x \leq \beta}  q_1\left(x\right) \leq C, \qquad \sup_{ \frac{1}{ \beta} \leq x \leq \beta} \int_1^ \infty ds q_s \left( x\right) \leq C.
\end{equation}
This is a consequence of the fact that $q_{1}$ is bounded for the first inequality, and the second inequality follows from the scaling property \eqref{eq:scaling} combined with the fact that $p_{1}$ is bounded (see e.g. \cite[I. 4.]{Zolotarev}). To establish \eqref{eq:approx2}, we first show that
\begin{equation}
\label{eq:affreux}\lim_{n \rightarrow \infty }
\left( \sup_{ \frac{1}{\a} B_n \leq k \leq \a B_n} \left|
\frac{\phi_{n- \fl{a n}}(k+1)/\phi_n(1)}{\phi^*_{n-\fl{a
n}}(k+1)/\phi^*_n(1)}- \theta \frac{ \frac{1}{1-a} q_1 \left( \frac{k+1}{B_ {n- \fl {an}}} \right)}{ \int_1^ \infty ds q_s \left( \frac{k+1}{ B_ {n- \fl {an}}}\right) } \right|\right)=0
\end{equation}
Since $B_{n- \fl{a n}}/B_{n} \rightarrow (1-a)^ {1/  \theta}$ as $n \rightarrow \infty$, Lemma \ref {lem:technical2} garanties the existence of two sequences $(\varepsilon^{(1)}_{k,n}, \varepsilon^{(2)}_{k,n})_{k, n \geq 1}$ such that
$$(n- \fl{a n}) \phi_{n- \fl{a n}}(k+1)=
q_1\left(\frac{k+1}{B_{n- \fl{a n}}}\right) +\varepsilon^{(1)}_{k,n},   \qquad \phi^*_{n- \fl{a n}}(k+1)= \int_1 ^\infty ds \, q_s\left(\frac{k+1}{
 B_{n- \fl{a n}} }\right) +\varepsilon^{(2)}_{k,n}$$
and such that $\max (\varepsilon^{(1)}_{k,n}, \varepsilon^{(2)}_{k,n}) \rightarrow0$ as $ n \rightarrow \infty$, uniformly in  ${1}/{\a} \cdot B_n \leq k \leq \a B_n$.  To simplify notation set $ m_{n}=n- \fl{an}$. By \eqref{eq:C} and the fact that $B_{m_{n}}/B_{n} \rightarrow (1-a)^ {1/  \theta}$, there exists $C>0$ such that for $n$ sufficiently large and ${1}/{\a} \cdot B_n \leq k \leq \a B_n$:
\begin{eqnarray*}
\left| \frac{m_n\phi_{m_n}(k+1)}{\phi^*_{m_n}(k+1)}-  \frac{  q_1 \left( \frac{k+1}{B_ {m_n}} \right)}{ \int_1^ \infty ds q_s \left( \frac{k+1}{ B_ {m_n}}\right) }  \right|
&=& \left|\frac{\varepsilon^{(1)}_{k,n} \cdot \int_1^ \infty ds q_s \left( \frac{k+1}{ B_ {m_n}}\right)- \varepsilon^{(2)}_{k,n} \cdot  q_1\left(\frac{k+1}{B_{m_{n}}}\right)}{ \int_1^ \infty ds q_s \left( \frac{k+1}{ B_ {m_n}}\right) \cdot \left(  \int_1^ \infty ds q_s \left( \frac{k+1}{ B_ {m_n}}\right)+ \varepsilon^{(2)}_{k,n}\right)} \right| \\
& \leq &  \frac{2C}{ \delta ^2} \cdot \sup_{ \frac{1}{\a} B_n \leq k \leq \a B_n} \max(\varepsilon^{(1)}_{k,n}, \varepsilon^{(2)}_{k,n}), \end{eqnarray*}
 where we have used \eqref{eq:delta} for the last inequality. This, combined with the fact that  $\phi^*_n(1)/(n\phi_n(1)) \rightarrow \theta$ as
$n \rightarrow \infty$ by Lemma \ref{lem:locallimit}, implies \eqref{eq:affreux}. Then our claim \eqref{eq:approx2} follows the scaling property \eqref{eq:scaling} and the continuity of $ \Gamma_{a}$.

We shall now prove another useful result before introducing some notation. Fix $\a>1$. Let $g_{n}: \R^ { \fl{an}+1} \rightarrow \R_{+}$ be a bounded measurable function. To simplify notation, for $x_{0}, \ldots, x_{\fl{an}}  \in \R$, set
$$G_{n}(x_{0}, \ldots, x_{\fl{an}})= g_{n}(x_{0}, \ldots, x_{\fl{an}}) \mathbbm{1}_{ x_{\fl{an}} \in \left[ \frac{1}{\a} B_n, \a
B_n \right]}$$
and, for a tree $ \tau$, set
$$ \widetilde{G}_{n}( \tau)=g_{n}( \W_0 ( \tau), \W_1 ( \tau), \ldots, \W_ { \fl {an}} ( \tau)) \mathbbm{1}_{  \{\W_{\fl{na}}(\tau) \in \left[ \frac{1}{\a} B_n, \a
B_n \right]\}}.$$   We claim that
\begin{equation}\label{eq:utileclaim}\lim_{n \rightarrow \infty} \left| \Es{  \widetilde{G}_{n}(\t_n) }-  \Esmu{ \left.  \widetilde{G}_{n}(\tau) \Gamma_a \left(\frac{ \W_{\fl{an}}(\tau)}{B_n} \right) \right| \, \zt \geq n} \right|=0.\end{equation}
Indeed, using successively Proposition \ref{prop:RW} and (\ref{eq:rel2}), we have:
\begin{eqnarray*}
 &&\Es{  \widetilde{G}_{n}(\t_n) }- \notag \Esmu{ \left.  \widetilde{G}_{n}(\tau) \Gamma_a \left(\frac{ \W_{\fl{an}}(\tau)}{B_n} \right) \right| \, \zt \geq n} \\
 && \qquad =	\Es{G_n(W_0, \ldots, W_ {\fl{na}}) |\, \z_1 =n}- \Es{G_n(W_0, \ldots, W_ {\fl{na}}) \Gamma_a \left(\frac{ W_{\fl{an}}}{B_n} \right) |\, \z_1 \geq n} \\
 && \qquad = \Es{ \left. G_n(W_0, \ldots, W_ {\fl{na}} ) \left( \frac{\phi_{n-[a
n]}(W_{\fl{an}}+1)/\phi_n(1)}{\phi^*_{n-\fl{an}}(W_{ \fl{a n}}+1)/\phi^*_n(1)} - \Gamma_a \left(\frac{ W_{\fl{an}}}{B_n} \right)\right)\right|\, \z_1 \geq n}.
\end{eqnarray*}
Our claim \eqref{eq:utileclaim} then follows from (\ref{eq:approx2}).

We finally return to the proof of Lemma \ref {lem:cvH0a}. Let $F : \D([0,a],\R) \rightarrow \R_+$ be a bounded
continuous function. We also set $F_n(\tau)=F\left(\frac{B_n}{n} H_{ \fl{n t}}(\tau);\, 0 \leq t \leq a \right)$. Since $(H_0 ( \tau), H_1 ( \tau), \ldots, H_ { \fl {an}} ( \tau))$ is a measurable function of $( \W_0 ( \tau), \W_1 ( \tau), \ldots, \W_ { \fl {an}} ( \tau))$ (see \cite[Prop 1.2]{RandomTrees}), by \eqref{eq:utileclaim} we get:
\begin{equation}\lim_{n \rightarrow \infty} \left| \Es{  F_n(\t_n) \mathbbm{1}_{ A_n^ \a ( \t_n)}}- \notag \Esmu{ \left.  F_n(\tau) \mathbbm{1}_{ A_n^ \a ( \tau)} \Gamma_a \left(\frac{ \W_{\fl{an}}(\tau)}{B_n} \right) \right| \, \zt \geq n}
\right|=0.\label{eq:cvutiliser2}\end{equation}
 By Theorem \ref{thm:cvW}, the law of
$$\left(\left(\frac{B_n}{n} H_{ \fl{nt}}(\tau);\, 0 \leq t \leq a
\right),  \frac{1}{B_n} \W_{ \fl{an}}(\tau) \right)$$ 
under $\Pmuzgeqn$  converges towards the law of $((H_t; 0 \leq t \leq a), X_a)$
under $\Nn(\cdot \| \z>1)$ (for the convergence of the second
component we have also used the fact that $X$ is almost surely continuous at
$a$). Thus:
\begin{eqnarray} \lim_{n \rightarrow \infty} \Es{F_n(\t_n) \mathbbm{1}_{
\left\{ \W_{\fl{na}}(\t_n) \in \left[ \frac{1}{\a} B_n, \a
B_n \right]
\right\}}}
&=&\Nn \left( F(
H_t; 0 \leq t \leq a) \Gamma_a(X_a) \mathbbm{1}_{\{ X_a \in \left[\frac{1}{\a},
\a \right] \} }\| \z
> 1\right) \notag\\
&=&\Nn \left( F(
H_t; 0 \leq t \leq a) \mathbbm{1}_{\{ X_a \in \left[\frac{1}{\a},
\a \right] \} }\| \z
= 1\right), \label{eq:utile}
\end{eqnarray}
where we have used  Proposition \ref{prop:Gamma} in the second equality.

By taking $F\equiv 1$, we obtain: $$\lim_{n \rightarrow \infty} \Pr{ \W_{\fl{na}}(\t_n) \in \left[ \frac{1}{\a} B_n, \a
B_n \right]} = \Nn\left( \left.X_a \in \left[\frac{1}{\a},
\a \right] \right| \, \z=1\right).$$
This last quantity tends to $1$ as $ \alpha \rightarrow \infty$.  By choosing $ \a>1$ sufficiently large, we easily deduce from the convergence \eqref{eq:utile} that:
$$ \lim_{n \rightarrow \infty} \Es{F\left(\frac{B_n}{n} H_{ \fl{n t}}(\t_n);\, 0 \leq t \leq a \right)} = \Nn \left( F(
H_t; 0 \leq t \leq a) \| \, \z
= 1\right). $$
The path continuity of $H$ under $\Nn \left(  \, \cdot \, \| \, \z
= 1\right)$ then implies the claim of Lemma \ref {lem:cvH0a}.
\end{proof}

\begin{theorem}Let $\t_n$ be a random tree distributed according to $\Pmuzn$. Then:
 $$\left(\frac{B_n}{n} H_{nt}(\t_n), \frac{B_n}{n} C_{2nt}(\t_n)\right)_{0 \leq t \leq 1} \qquad
 \mathop{\longrightarrow}_{n \rightarrow \infty}^{(d)} \qquad
 (H_t,H_t)_{0 \leq t \leq 1}  \textrm{ under } \Nn( \, \cdot\, | \z
 =1).$$\end{theorem}

\begin{proof}
The proof consists in showing that the scaled height process is
close to the scaled contour process and then using a
time-reversal argument in order to show that the convergence holds
on the whole segment $[0,1]$. To this end, we adapt \cite[Remark
3.2]{Duquesne} and \cite[Section 2.4]{DuquesneLG} to our context.
For $0 \leq p < {n}$ set $b_p=2p-H_p(\t_n)$ so
that $b_p$ represents the time needed by the contour process to
reach the $(p+1)$-st individual of ${\t_n}$ (in the lexicographical
order). Also set $b_{n}=2 ({n}-1)$. Note that $C_{b_p}=H_p$. From this observation, we get:
\begin{equation}\label{eq:ineg}\sup_{t \in [b_p,b_{p+1}]}
|C_t(\t_n)-H_p(\t_n)| \leq
|H_{p+1}(\t_n)-H_p(\t_n)|+1.\end{equation} for $0 \leq p < n$. Then
define the random function $g_n : [0, 2 {n}] \rightarrow \N$ by
setting $g_n(t)=k$ if $t \in [b_k,b_{k+1})$ and $k<{n}$, and
$g_n(t)={n}$ if $t \in [2({n}-1),2 {n}]$ so that for $t < 2(n-1)$, $g_n(t)$ is the index
of the last individual which has
been visited by the contour function up to time $t$ if the individuals are indexed $0,1, \ldots, n-1$
in lexicographical order. Finally, set
$\widetilde{g}_n(t)=g_n({n} t)/{n}$. Fix $ a \in (0,1)$. Then,
by (\ref{eq:ineg}):
\begin{equation}
\label{eq:maj}\sup_{t \leq \frac{b_{\lfloor a n \rfloor}}{{n}}}\left| \frac{B_{n}}{{n}} C_{{n} t}(\t_n)- \frac{B_{n}}{{n}} H_{{n} \widetilde{g}_n(t)}( \t_n)\right|
\leq \frac{B_{n}}{{n}}+\frac{B_{n}}{{n}} \sup_{k \leq \lfloor a
n \rfloor}|H_{k+1}(\t_n)-H_k(\t_n)|,
\end{equation} which converges in
probability to $0$ by Lemma \ref{lem:cvH0a} and the path continuity of $(H_t)$. On the other hand it follows from
the definition of $b_n$ that: $$\sup_{t \leq \frac{b_{\lfloor
a n \rfloor}}{{n}}}\left| \widetilde{g}_n(t)-\frac{t}{2}
\right| \leq \frac{1}{2 B_{n}} \sup_{k \leq a {n}}
\frac{B_{n}}{{n}}
H_k(\t_n)+\frac{1}{{n}}\quad\mathop{\longrightarrow}^{(\P)}\quad 0$$
 by Lemma \ref{lem:cvH0a}. Finally, by the definition of $b_n$ and using
 Lemma \ref{lem:cvH0a} we see that $\frac{b_{\lfloor a n
\rfloor}}{{n}}$ converges in probability towards $2a$ and that $ \frac{B_n}{n} \sup_ {t \leq 2a} \left| H_ {{n} \widetilde{g}_n(t)}( \t_n)-H_ {nt/2}( \t_n)\right|$ converges in probability towards $0$ as $n \rightarrow \infty$. Using \eqref{eq:maj}, we conclude
that:
\begin{equation}\label{eq:HC}\frac{B_{n}}{{n}} \sup_{0 \leq t \leq a} |C_{2 {n} t}(\t_n) - H_{{n} t}(\t_n)|
\quad \mathop{\longrightarrow}^{(\P)} \quad 0.\end{equation}
Together with  Lemma \ref{lem:cvH0a}, this implies:
$$\left(\frac{B_{n}}{{n}} C_{ 2 n
t}(\t_n);\, 0 \leq t \leq a \right)  \quad
\mathop{\longrightarrow}^{(d)} \quad  \left( H_t;\, 0 \leq t \leq a
\right) \textrm{ under }\Nn( \, \cdot \, | \, \z=1).$$ Since $ (C_t(\t_n); \, 0 \leq t \leq 2n-2)$ and
$(C_{2n-2-t}(\t_n); \, 0 \leq t \leq 2n-2)$ have the same
distribution, it follows that:
\begin{equation}\label{eq:cvC}\left(\frac{B_{n}}{{n}} C_{ 2 n
t}(\t_n);\, 0 \leq t \leq 1 \right)  \quad
\mathop{\longrightarrow}^{(d)}  \quad \left( H_t;\, 0 \leq t \leq 1
\right) \textrm{ under }\Nn( \, \cdot \, | \, \z=1).\end{equation}
See the last paragraph of the proof of Theorem 6.1 in \cite{LGIto}
for details.

Finally, we show that this convergence in turn entails the
convergence of the rescaled height function of $\t_n$ on the whole
segment $[0,1]$. To this end, we verify that convergence (\ref{eq:HC}) remains
valid for $a=1$. First note that:
\begin{equation}
\label{eq:cvg}\sup_{0 \leq t \leq 2 }\left| \widetilde{g}_n(t)-\frac{t}{2}
\right| \leq \frac{1}{2 {n}} \sup_{k \leq {n}}
H_k(\t_n)+\frac{1}{{n}}= \frac{1}{2 B_{n}} \sup_{k \leq 2 {n}}
\frac{B_{n}}{n} C_k(\t_n)+\frac{1}{{n}}
\quad\mathop{\longrightarrow}^{(\P)}\quad 0
\end{equation}
 by (\ref{eq:cvC}).
Secondly, it follows from \eqref{eq:ineg} that: \begin{eqnarray*} \sup_{0 \leq t \leq 2}\left|
\frac{B_{n}}{n} C_{{n} t}(\t_n)- \frac{B_{n}}{n} H_{{n}
\widetilde{g}_n(t)}\right| &\leq& \frac{B_{n}}{n} + \frac{B_{n}}{n}
\sup_{k \leq n-1
}|H_{k+1}(\t_n)-H_k(\t_n)| \\
&=& \frac{B_{n}}{n}+\frac{B_{n}}{n} \sup_{k \leq n -1}
\left|C_{\frac{b_{k+1}}{n}n}(\t_n)-C_{\frac{b_{k}}{n}n}(\t_n)\right|.\end{eqnarray*}
By (\ref{eq:cvC}), in order to prove that the latter quantity tends
to $0$ in probability, it is sufficient to verify that $\sup_{k \leq n
} \left| \frac{b_{k+1}}{n}-\frac{b_{k}}{n}\right|$ converges to $0$
in probability. But by the definition of $b_n$:
\begin{eqnarray*}\sup_{k \leq n } \left|
\frac{b_{k+1}}{n}-\frac{b_{k}}{n}\right|=\sup_{k \leq n } \left|
\frac{2+H_k(\t_n)-H_{k+1}(\t_n)}{n}\right| &\leq&
\frac{2}{n}+ 2\sup_{k \leq n } \frac{H_k(\t_n)}{n}
\end{eqnarray*}
which converges in probability to $0$ as in \eqref{eq:cvg}.
As a
consequence:
$$\frac{B_{n}}{n} \sup_{0 \leq t \leq 1} |C_{2 {n} t}(\t_n) - H_{{n} \widetilde{g}_n(2t)}(\t_n)|
\quad \mathop{\longrightarrow}^{(\P)} \quad 0.$$
By (\ref{eq:cvC}), we get that:
$$\left(\frac{B_n}{n}H_{{n} \widetilde{g}_n(2t)}(\t_n)\right)_{0 \leq t \leq 1 } \quad\mathop{\longrightarrow}^{(d)}_{n \rightarrow \infty}\quad
(H_t)_{0 \leq t \leq 1} \textrm{ under }\Nn( \, \cdot \, | \,
\z=1).$$
Combining this with
(\ref{eq:cvg}), we conclude that:
$$\left(\frac{B_n}{n}C_{2 n t}(\t_n), \frac{B_n}{n}H_{n t}(\t_n)\right)_{0 \leq t \leq 1 } \quad
\mathop{\longrightarrow}^{(d)}_{n \rightarrow \infty}\quad (
H_t,H_t)_{0 \leq t \leq 1} \textrm{ under }\Nn( \, \cdot \, | \,
\z=1).$$
This completes the proof.
\end{proof}

\begin{remark}If we see the tree $\t_n$  as
a finite metric space using its graph distance, this theorem implies that  $\t_n$,  suitably rescaled,
converges in distribution to the $\theta$-stable tree, in the sense of the
Gromov-Hausdorff distance on isometry classes of compact metric
spaces (see e.g. \cite[Section 2]{RandomTrees} for details).
\end{remark}

\begin{remark}When the mean value of $\mu$ is greater than one, it is possible to replace $\mu$
with a critical probability distribution belonging to the same
exponential family as $\mu$ without changing the distribution of
$\t_n$ (see \cite{Kennedy}). Consequently, the theorem holds in the
supercritical case as well. The case where $ \mu$ is subcritical and $ \mu(i) \sim L(i)/i^ { 1+\theta}$ as $i \rightarrow \infty$ has been treated in \cite {Knongeneric}. However, in full generality, the non-critical subcritical case remains open.
\end{remark}

\subsection{Proof of the technical lemma}

In this section, we prove Lemma \ref{lem:technical2}.

\begin {proof}[of Lemma \ref{lem:technical2}]
We first prove (i). By the local limit theorem (Theorem \ref{thm:locallimit} (ii)),
we have, for $k \geq 1$ and $j \in \Z$:
$$  \left|{B}_n \P[W_n=j]-
p_1\left(\frac{j}{{B}_n}\right) \right| \leq\epsilon(n),$$
where $ \epsilon(n) \rightarrow 0$. By Proposition \ref{prop:kemperman}, we have
$
 n \phi_n(j) = j \P[W_n=-j]$. Since $ \frac{ j }{{B}_n} p_1\left(-\frac{j}{{B}_n}\right)=q_1\left(\frac{j}{B_{n}}\right)$, we have for $1 \leq j \leq \a B_n$:
 $$\left|n
 \phi_n(j)-q_1\left(\frac{j}{B_{n}}\right) \right|= \frac{j}{{B}_n} \left|  B_n \P[W_n=-j]-  p_1\left(\frac{j}{{B}_n}\right)\right|\leq \a \e(n).$$
 This
 completes the proof of (i).

For (ii), first note that by the definition of $q_s$ and the scaling property (\ref{eq:scaling}):
$$\int_1 ^\infty ds \, q_s\left(\frac{j}{B_n}\right)= \int_1 ^ \infty  \frac{j/B_n}{s^{1/\theta+1}} p_1\left(- \frac{j/B_n}{s^{1/\theta}}\right)
\,ds.$$
 By Proposition
\ref{prop:kemperman} and the local limit theorem:
$$
\left|\phi^*_n(j) - \sum_{k=n}^\infty  \frac{ j }{ k {B}_k}
p_1\left(-\frac{j}{{B}_k}\right) \right|= \left| \sum_{k=n}^\infty \left( \frac{j}{k} \P[W_k=-j] - \frac{ j }{ k {B}_k}
p_1\left(-\frac{j}{{B}_k}\right) \right) \right| \leq \sum_{k=n}^\infty \frac{ j
}{k {B}_k}
 \e(k),$$
 where $\e(n)\rightarrow 0$.
Then write:
\begin{eqnarray*}&& \left|\sum_{k=n}^\infty  \frac{ j }{ k {B}_k}
p_1\left(-\frac{j}{{B}_k}\right)-\int_1 ^ \infty ds \,
\frac{j/B_n}{s^{1/\theta+1}} p_1\left(-
\frac{j/B_n}{s^{1/\theta}}\right)  \right|\\
&& \quad \leq \int_1 ^ \infty ds \, \left| \frac{j n}{B_{\fl{ns}}
\fl{ns}} - \frac{j/B_n}{s^{1/\theta+1}} \right| p_1\left(-
\frac{j}{B_{\fl{ns}}}\right) +\int_1 ^ \infty ds \,
\frac{j/B_n}{s^{1/\theta+1}}\left| p_1\left(-
\frac{j}{B_{\fl{ns}}}\right)- p_1\left(-
\frac{j/B_n}{s^{1/\theta}}\right)\right|.
\end{eqnarray*}
Denote the first term of the right-hand side by $P(n,j)$ and the second term
by
$Q(n,j)$ . Since $p_1$ is bounded by a constant
which we will denote by $M$, we have for $1 \leq j \leq \a B_n$:
\begin{eqnarray*}
 P(n,j) \leq \a M \int_1 ^ \infty ds \, \frac{1}{s^{1/\theta+1}}\left|
\frac{n B_n s^{1/\theta+1}}{B_{\fl{ns}} \fl{ns}} -1 \right|.
\end{eqnarray*}
For fixed $s \geq 1$,$\frac{1}{s^{1/\theta+1}}\left|
\frac{n B_n s^{1/\theta+1}}{B_{\fl{ns}} \fl{ns}} -1 \right|$ tends to $0$ as $n \rightarrow \infty$, and using Proposition \ref{prop:slow}, the same quantity is bounded by an integrable function independent of $n$.
The dominated convergence theorem thus shows that $ P(n,j)
\rightarrow 0$ uniformly in $1 \leq j \leq \a B_n$.
Let us now bound $Q(n,j)$ for  $1 \leq j \leq \a B_n$. Since the
absolute value of the derivative of $p_1$ is bounded by a
constant which we will denote by $M'$, we have:
$$Q(n,j) \leq M' \int_1 ^ \infty ds \,
\frac{j/B_n}{s^{1/\theta+1}}\left|
\frac{j}{B_{\fl{ns}}} -
\frac{j/B_n}{s^{1/\theta}}\right| \leq\a^2 M' \int_1 ^ \infty ds \, \frac{1}{s^{2/\theta+1}}
\left| \frac{B_n s^{1/\theta}}{B_{\fl{ns}}}-1\right|.$$ The right-hand side tends to $0$ by the same
argument we used for $P(n,j)$. We have thus proved that:
$$ \lim_{n \rightarrow \infty}\sup_{1 \leq j \leq \a B_n} \left|\sum_{k=n}^\infty  \frac{ j }{ k {B}_k}
p_1\left(-\frac{k}{{B}_k}\right)- \int_1 ^\infty ds \,
q_s\left(\frac{j}{
 B_{n} }\right)\right|=0.$$
One finally shows that $ \sum_{k=n}^\infty \frac{ j}{k {B}_k} \e(k)$ tends to $0$ as $n \rightarrow \infty$ uniformly in $1 \leq j \leq \a B_n$
 by noticing that:
 $$  \sup_ {n \geq 1} \sup_{1 \leq j \leq \a B_n} \left( \sum_{k=n}^\infty \frac{ j
}{k {B}_k}
\right) \leq \a \sup_ {n \geq 1}\left( \sum_{k=n}^\infty \frac{ B_n
}{k {B}_k}
\right)  <\infty.$$
 This completes the proof.
\end {proof}

\begin {tabular}{l }
Laboratoire de mathmatiques, \\
UMR 8628 CNRS.\\
Universit Paris-Sud\\
91405 ORSAY Cedex, France
\end {tabular}

\medbreak
\noindent \texttt{igor.kortchemski@normalesup.org}
\end{document}